\documentclass[12pt]{amsart}
\usepackage{latexsym}
\usepackage{amsfonts}
\usepackage{amsthm}
\usepackage{amsmath}
\usepackage{amssymb}

\def\<{{\langle}}
\def\>{{\rangle}}

\def\note#1{{}}

\def\note#1{}

\def\beq{\begin{equation}}
\def\eeq{\end{equation}}

\newcounter{zlist}

\def\Label#1{\label{#1}\ifmmode\llap{[#1] }\else
\marginpar{\smash{\hbox{\tiny [#1]}}}\fi}
\def\Label{\label}

\newtheorem{proposition}{Proposition}[section]

\newtheorem{theorem}[proposition]{Theorem}

\theoremstyle{definition}
\newtheorem{definition}[proposition]{Definition}

\theoremstyle{remark}
\newtheorem{remark}[proposition]{Remark}

\newcounter{c}

\newcommand{\etyk}[1]{\vspace{-7.4mm}$$\begin{equation}\Label{#1}
\addtocounter{c}{1}}
\renewcommand{\]}{\ifnum \value{c}=1 $$\else \end{equation}\fi}
\begin{document}

\title{Remembering Solomon Marcus}
\author{Florin Felix NICHITA}
\address{Institute of Mathematics of the Romanian Academy, 
P.O. Box 1-764, RO-70700, Bucharest, Romania}
\date{}
\begin{abstract} {From the Andre Breton's manifest
of surrealism,  through the transdisciplinary understanding,
we arrive at
``A post-modern manifest'' (Sci 2022).
A talk  by
Laura De Marco (Harvard)  will provide scientific background
to approach an
AMS poetry.
The next section  will be a qualitative 
analysis of some new
operations on the real numbers.
 The conclusions
will be given in the last section, and an appendix
will recall some previous work with some new comments.
 }
\end{abstract}  
\maketitle

\section{Introduction}

Solomon Marcus extended his activities to many domains.
Meeting my university colleagues and me,
 he remarked the Ph.D. Thesis of Nicusor Dan,
the achievements of Marian Aprodu,
 the contributions
of Gabriel Istrate in Computer Science
(although he graduated as a pure mathematician), 
the carrier of Ovidiu Calin (a Fulbright Fellow in Romania),
the Euler type formulas found by myself,
etc.

Our  papers on meetings with Solomon Marcus, 
which might become  a book,
received excellent feedback from many readers. We present
an aftermath of that series.
 Trying to recreate
the spirit of those meetings, we 
involuntarily gathered an entire ``choir'' on many ``voices'':
Horatiu Malaele,
Mircea Cartarescu, Sanda Golopentia, Rodica Zafiu,
Gabriel Prajitura (the literature voices),
 Gheorghe Gussi, Ovidiu Pasarescu,  Marius Buliga, Remus Radu,
Ovidiu Calin, Radu Iordanescu (the math choir),
Violeta Dinescu, Nick Mihalache, Dan Vuza, Aurelian Gheondea,
Ion Mihail Nichita, Iosif Imbrisca
(the ``music department'')
and the Math Cafe Team (Naomi,
David, Razvan and Rarea, among others).
From Marius Buliga, I have learned that the
mathematics archive (``arXiv.org'') posts 
many papers of Leonard Euler in English (in relation with
the paper \cite{SM}).
Also, I found at ``Project Euler''  the following sum:
$ \sum^{10}_1 k^k $  (this would be related to
 \cite{mMemo2}).

Just as the Andre Breton's manifest of surrealism
led to masterpieces of artwork combining
the real life with the dreams, 
in mathematics,  the real numbers
and the imaginary numbers live in the field of complex numbers.
Combining the arithmetic with the world of imaginary numbers,
 the Euler identity,
$ \  \   e^{i \pi} + 1 = 0 \ , \  $
 is an example
of ``surrealistic'' masterpiece in mathematics.
The Euler formula and the Euler identity were
our favourite topics in our meetings with Solomon Marcus.
These led to a publication in Axioms (\cite{SM}).

The
``Manifesto of Transdisciplinarity'' of Basarab Nicolescu,
encouraged me  for the involvement  in many disciplines.
Solomon Marcus could be considered  an example
in this regard.
The poem
``A post-modern manifest'' (Sci, 2022)  
invites the poets to  explore the mathematical universe.
Now, I propose a 
 short intermezzo in our discourse: 
Laura De Marco, from Harvard University, gave a talk on
the Mandelbrot set  at the 
Romanian Academy in June 2026. 
This vibrant presentation attracted the attention of
my colleagues. Four collaborators of Laura De Marco
are Romanians. And a final comment for this intermezzo:
Mandelbrot was a polymath.
Back to  our main discourse.  The
AMS poetry 
 ``Pantoum for the Mandelbrot set'' 
won the second prize in 2026. It
is both
full of lyricism and of scientific information.

The next part of the paper is an 
analysis of 
operations defined on the set of real numbers.
The conclusions
will be given in the last section, and an appendix
will recall some previous work and some possible continuations.


\section{Mathematical background}


\subsection{New logarithms and radicals}
We will define an ``universal logarithm'', a kind of logarithm 
for which there is no need to give
a fixed basis.
Thus, we denote by
$ L_2 $ the universal logaritm of order two:  given a number $n$,
the universal logarithm associates the bigest number $m$, such that
$ n = m^m $.
For example:
$  L_2 (4) = 2 $ and
$  L_2 (27) = 3 $.

The universal
logarith of order three has the following values:\\
$ L_3 (16) \ = L_3 (2^{2^2}) \ = 2 $,
$ \ L_3 (1) \ = L_3 (1^{1^1}) \ = 2 $,
$  \ L_3 (10^{10^{10}}) \ = 10 $,


We also define a generalized radical: it takes
the radical of the exponent.

$ R_{{2}} \ (x) =  \sqrt{2}^{\sqrt{\log_{\sqrt{2}} x}} \ \ . \ \ $
For example, $ R_{{2}} \ (4) =  \sqrt{2}^{\sqrt{\log_{\sqrt{2}} 4}}
=  \sqrt{2}^{\sqrt{4}} =  \sqrt{2}^{2} = 2 $.

In general, $  \ \ L_2 (x) \ \leq \  R_{{2}} \ (x) $.


\subsection{New operations}

The following operation could be considered
an operation which (approximately) generates the addition of
 positive real numbers:

$  \ \ \ \ \ \ \ 
a \circ b = \log_{\sqrt{2}} \ {\sqrt{2}}^{ \ a}  \ + 
\ {\sqrt{2}}^{ \ b}
 \ . \ \ $
For example, $ \ \ \ \ a \circ a \circ a \circ  a = a + 4 $.


The following operation also generates the addition of
 positive real numbers:\\
$ \ \ \ \ \ \ \ \ \ \ \ a  \oplus a = a+2 \ ; $
$  \ \ \ 
a \oplus b = \max(a, b) + 1 \ , \ a \neq b  \ . \ \ $
For example, $ \ \ \ \ a \circ a \circ a  = a + 3 $.



We now define  a symmetric exponentiation:

$ \ \widehat{a^b} = \ \widehat{b^a} =
\sqrt{2}^{( \log_{\sqrt{2}} a) \times ( \log_{\sqrt{2}} b)} =  
   a^{ \log_{\sqrt{2}} b} = b^{ \log_{\sqrt{2}} a} $.

For example, $ \ \widehat{  2^2} = 4  \ $,
$ \ \widehat{  2^{\sqrt{2}}} = 2 \ \ $,
$ \ \widehat{  2^{1}} = 1 \ \ $
 and 
 $ \ \widehat{  2^4}  = 16 $.

\subsection{Mean type inequalities}
For $ \ \ a, b \geq 4 $, we have the following inequalities:

$$ L_2 ( a^b) \leq \sqrt{ab} \leq \frac{a+b}{2} \leq a \circ b - 
2 $$


For $ \ \ a, b \geq 1 $, natural numbers, we have the following inequalities:

$$ R_{{2}} \ ( \widehat{a^b}) \  \leq \sqrt{ab} \leq \frac{a+b}{2} \leq a \oplus
b - 
\frac{3}{2} $$


\subsection{Euler formula}
There exists
 an Euler formula for the set of real numbers endowed with
 the operations $ \ \oplus $ and $ \ + $ . A new element
$ \ - \infty $ will be added to the real numbers, and it is
a neutral elemnent for $ \ \oplus $. 

The ``exponential function'' $ f(t) = \sqrt{2}^{  \ t } $, 
$ \  \sqrt{2}^{ \ s \oplus t } = \sqrt{2}^{ \ s } + \sqrt{2}^{  \ t } $, 
 can be extended to 
$ \mathbb{R} \times \mathbb{R}$ as follows:
 $ \ \ \ \  \sqrt{2}^{ \ (- \infty,t)} \ =\ [ \ \log_2 \ ( \cos^2 \sqrt{2}^{ \ t } ) \ , 
\ \log_2 \ ( \sin^2 \sqrt{2}^{ \ t } ) \ ] \ . $

\bigskip

\section{Commentaries and conclusions}


Some traces of surrealism could be found in  (the pictures of) 
the famous books
of Lewis Carroll. Forms of contemporary surrealism are present
in 
 the magical 
realism (of Gabriel Garcia Marquez, Jorge Luis Borges, etc)
and in the works of the street artist Banksy (probably
an acronym for the ``Bristol
underground scene''). A dual  ``art'' has emerged: finding
Banksy's works across the globe, admiring them,
capturing them,
transporting them
 and exposing them in official  exhibitions.
(The last three steps are not fully accepted by Banksy and
his admires.)
Other contemporary artists created Dog Man, 
present the aboriginal culture
or use special
organic pigments in their ``surrealistic art''
(Bojaxhiu, Yuguo, etc).

\bigskip
\bigskip

\section{APPENDIX -- Previous work revisited}

I received an invitation from the journal
``Libertas Matematica 
(new series)'' to write an article
about Solomon Marcus 
on the occasion of one hundred years
from his birth.
I
submitted my
 manuscript (\cite{Memo}), and it
 was 
accepted soon afterwards.
This determined me to write a second paper
dedicated to Solomon Marcus (\cite{mMemo}) . 
Also,
 I posted a survey of  \cite{Memo} and \cite{mMemo} 
on the 
Mathematical Archive (\cite{mMemo2}), adding
some new content. (For example, I
gave applications of the B-ring Euler formula in
finding solutions for the braid condition.)
From Brown University (Rhode Island), I received four
articles about Solomon Marcus (see \cite{s1,s2,s3,s4}).
Also, I would like to
mention another Mathematical Archive communication (\cite{arxive}),
which was posted soon after our preprint (\cite{mMemo2}), and it gives some
weight
to our results, because the solutions to the braid equation
lead to representations of the braid group. More recently,
Bradshaw and Vignat wrote a paper, in a similar manner
with ours, about another
``beautiful mind'' 
(see \cite{bv}). Also, Terence Tao has posted an essay on 
Mathematics and AI, 
 reminding me of a book by Marta Petru 
and a series of articles by Calin Vlasie
(see \cite{tt, mp, cv}). The International Poetry Festival in 
Bucharest (FIPB 2026) is a sourse of wonderful
discussions on the topics of our paper.

\subsection{Trigonometry   in racks}

\begin{definition}
 A rack is triple $ ( S, \ \cdot, \ \diamond \ )$, where
$ S $ is
a set with two binary operations,
satisfying the following  four axioms:

$ a \cdot (b \cdot c) = (a \cdot b) \cdot (a \cdot c) \ , $
$  \ \ \  \ \ \  \ \ \ \ \ \ \ \ \ \ \ \ (a \cdot b) \  \diamond  \  a = b \ ,$

$ 
a \cdot (b \ \diamond  \  a)  = b \ $,
$ \ \ \ \ \ \ \ \ \ \ 
 ( c \ \diamond \  b) \ \diamond  \ a = 
( c \ \diamond  \  a) \  \diamond  \  ( b \ \diamond  \  a) \ .$

The operation $ \  \cdot \ $ is called the main operation, we will
write $ a \cdot b = ab \ $, and it has priority over $ \ \diamond  \ $
 in formulas.
\end{definition}


\begin{remark}
 The following is an example of rack associated to a group:
 
$ \ \ ab = aba^{-1} $, $ \ \ a \ \diamond  \   b = b^{-1}ab $.
\end{remark}

\bigskip

We now choose $ \ \ \ \ e, \ O  \in S$, and let
$ \ \ \Pi = eO \ \ $ and 
$  \ U = e(eO) $.

Let $ \cos x = ex \ , \ \sin x = x \ \diamond  \  e \ \ $ be 
``trigonometric'' functions
in our rack.


The following properties hold for our ``trigonometric'' functions:

$ \ \ \ \ \ \cos \Pi = U \ ; \ \ \ \ \ \ \ \sin \Pi = O \ $;


$ \cos \ xy =  \cos x \ \cos y \ ;
\ \ \ \ \ \ \ \ \ \ \ \ \ \ 
 \cos (x \ \diamond  \ y) =  \cos x \ \diamond  \  \cos y \ $;

$ \sin \ x   y = \sin x \ \sin y \ ;
\ \ \ \ \ \ \ \ \ \ \ \ \ \ 
\sin (x \ \diamond  \ y) = \sin x \ \diamond  \  \sin y \ $.

The fundamental formula for this ``trigonometry'' is the following:
$$  \sin \ ( \cos x ) \ = \cos ( \sin x  ) = x. $$ 

\bigskip

\begin{remark}
 The following rack can be defined now:
$ \ \ ab = \cos b $, $ \ \ a \diamond   b = \sin a $.
\end{remark}

\subsection{Weak racks}

\begin{definition}
 A { \bf weak rack} is triple $ ( S, \ \cdot, \ \diamond \ )$, where
$ S $ is
a set with two binary operations,
satisfying the following  three axioms:

$ \ \  \ a \cdot (b \cdot c) = (a \cdot b) \cdot (a \cdot c) \ , $

$  
\ \  (a \cdot b) \  \diamond  \  a = a \cdot (b \ \diamond  \  a) \ ,$

$ \ \ 
 ( c \ \diamond \  b) \ \diamond  \ a = 
( c \ \diamond  \  a) \  \diamond  \  ( b \ \diamond  \  a) \ .$

Again, the for the main operation 
 we will
write $ a \cdot b = ab \ $, and it has priority over $ \ \diamond  \ $
 in formulas.
\end{definition}

\begin{remark}
 The following weak rack can be associated to a Boolean algebra:
$$ \ \ ab = a \rightarrow b ,  \ \ a \diamond   b =  a \setminus b \ . $$

Also, the following weak rack can be associated:
$ \ \ ab = a \vee b ,  \ a \diamond   b =  a \wedge b  . $
\end{remark}

\bigskip

Let $  \ e, \ O  \in S$, 
$  \ \Pi = eO  $, 
$  \ U = e(eO) $,
$ \cos x = ex \ ,$ and $ \ \sin x = x \ \diamond  \  e $. We have
the following properties:
$ \ \ \ \ \ \cos \Pi = U \ ; \ \ \ \ \ \ \ \sin \Pi = O \ $;
$ \cos \ xy =  \cos x \ \cos y \ ; $

$ \cos (x \ \diamond  \ y) =  \cos x \ \diamond  \  \cos y \ $;
$ \ \ \sin \ x   y = \sin x \ \sin y \ ;
\ \ \ \ \ \ \ \ \ \ \ \ \ \ 
\sin (x \ \diamond  \ y) = \sin x \ \diamond  \  \sin y \ $.

The fundamental formula is the following:
$  \ \ \sin \ ( \cos x ) \ = \cos ( \sin x  ) . $

\subsection{The Euler formula and the Euler identity   in 
(weak) racks}

\begin{definition} We can define a dual rack
with the opposite operations,
$ ( S, \ * \ ,  \bullet \ )$,
where the new operations are the following:
$ \ \ \ \  a * b = b \ \diamond \  a \ $ and 
$ \ \ a \bullet b = b \cdot a \ $.

\end{definition}

\begin{remark}
 The rack having the following operations is self-dual:
$ \ \ ab = b $, $ \ \ a \diamond   b = a $.
\end{remark}

On the set $ \ S \times S $, we can put a rack structure
obtained
as the  product of the initial rack with its dual.
Let us denote its main operation as follows:\\
$  \ \ {(x, \ y) \ \square \ (u, \ v) = (xu, \ v \diamond y)
= (xu, \  y * v)} $.


 For $ a \in S$, we  define 
an ``exponential'' function on $ \ S \times S $:\\
$ \ \  \exp_{a} (x, \ y) =
a^{(x, \ y)}  = (ax, \ a * y) $.



The ``exponential'' map has the property:
$ \ a^{(x, \ y) \square   (u, \ v)} = 
\ a^{(x, \ y)} \ \square \  a^{(u, \ v)}
 $.

The diagonal map, $ \ \Delta : S \rightarrow S \times S, \
x \mapsto ( x, \ x) \ $,  is a rack morphism for a certain rack structure
on $ \ S \times S$.

\bigskip

\begin{theorem} ({\bf Euler formula in racks.})
The following formula holds:
\begin{equation} 
\exp_{e} \circ \ \Delta = 
[ \cos  \ \times \ \sin ] \circ \Delta \ .
\end{equation}
Equivalently,
\begin{equation} 
{e}^{(x,\ x)} = 
( \cos x ,  \  \sin x )  .
\end{equation}

Moreover, the following identity is true:
$ \ \ \ \ \ \  e^{\Delta (\Pi )} = ( U,\ O ) \ . $

\end{theorem}

{ \bf Proof.} The left hand side reads:
$ \ \ \exp_{e} \circ \ \Delta (x) = e^{(x, \ x)} =
(ex, \ e*x) \ . $
The right hand side reads:
$ \ \ [ \cos  \ \times \ \sin ] \circ \Delta (x) 
= (\cos x, \ \sin x) \ . $
So, the left hand side equals the right hand side, because
$ \ \cos x = ex $ from the definition, and 
$ \ \sin x = x \diamond e = e*x \ $.

For the last identity we have:
$  e^{\Delta (\Pi )} =  e^{(\Pi , \ \Pi )} 
= ( e \Pi, \ e* \Pi) =  ( U,\ O )$.
\qed

\begin{remark}
The classical Euler's formula states that:
\begin{equation} \label{f}
  e^{ ix } = \cos x + i \sin x \ \ \ \ \ \forall x \in \mathbb{R}.
\end{equation} 

This formula could be also written as:
\begin{equation} 
\exp_{e} \circ \ j = 
[ \cos  \ \times \ \sin ] \circ \Delta \ ,
\end{equation}
where $ \ \ j : \mathbb{R} \rightarrow \mathbb{R} \times \mathbb{R} \ ,
\ \ x \mapsto (0, \ x) \ . $
For $x= \pi$ the Euler's formula becomes the Euler's identity. 
\end{remark}

\begin{remark} The conics are generated by:
\begin{equation} \label{f1}
  e^{ ix } = \cos x + i \sin x \ \ \ \ \ \forall x \in \mathbb{R}
\end{equation} 
\begin{equation} \label{f2}
  e^{ jx } = \cosh x + j \sinh x \  , \  j^2 =1,  \ \ \forall x \in \mathbb{R}
\end{equation} 
\begin{equation} \label{f3}
   x e^{ hx } = x ( 1 + h  x ) \  , \  h^2 =0,  \ \ \forall x \in \mathbb{R}.
\end{equation} 
\end{remark}

\begin{remark} The quadrics are generated by:
\begin{equation} \label{f11}
  e^{ x ( i \cos y + I \sin y) } = \cos x + i \sin x \cos y +
I \sin x \sin y \ , \
I^2 = -1 \ , \ iI=0 \ , \  \forall x \in \mathbb{R}
\end{equation} 
\begin{equation} \label{f11b}
  e^{ x ( i \cosh y + j \sinh y) } = \cos x + i \sin x \cosh y +
I \sin x \sinh y
\end{equation}
\begin{equation} \label{f111}
  e^{ x ( j \cos y + J \sin y) }  , \
J^2 = 1 \ , \ jJ=0 
\end{equation} 
\begin{equation} \label{f111b}
  x e^{  h \cos y + H \sin y) }  , \
H^2 = 0 \ , hH=0 
\end{equation} 
\begin{equation} \label{f1111}
  x^2 e^{  h \frac{\cos y}{x}  + H \frac{\sin y}{x} }  \ \ ,
\ \ x^2 e^{  h \frac{\cosh y}{x}  + H \frac{\sinh y}{x} }   \ .
\end{equation} 
\end{remark}

\bigskip

We  consider the following ``hyperbolic'' functions:\\
$ \ \cosh (x, y) = (ex,  y) \ $ and
$ \ \ \sinh (x, \ y) = (x, e * y) $.


It is easy to check the hyperbolic functions Euler formula:
$$ \  \ \exp_e  \ = \ \cosh \ \circ \ \sinh \ = 
\ \sinh  \  \circ \ \cosh  \ . $$

\bigskip

\begin{theorem}
The functions $ \ \exp_e \ , \ \cosh $ and $ \ \sinh $ are
solutions for the Quantum Yang-Baxter equation 
($ R^{12} \circ  R^{13} \circ R^{23} \ = \  
R^{23} \circ  R^{13} \circ R^{12} \  $). 
\end{theorem}


\bigskip

\begin{thebibliography}{Bibliography}{}

\bibitem{SM} Solomon Marcus and Florin Nichita, {\em On transcendental 
numbers: new results and a little history}, 
{ Axioms}, 2018, 7, 15.

 \bibitem{mMemo2} Nichita, F.,
{\em Memories with Solomon Marcus (III)}, arXiv: 2604.15970.



 \bibitem{Memo} Nichita, F.,
{\em Memories with Solomon Marcus}, Libertas Matematica 
(new series), Volume 45 (2025), No. 1, 103-108.

 \bibitem{mMemo} Nichita, F.,
{\em Memories with Solomon Marcus (II)}, in press.


 \bibitem{s1} Golopentia, S. {\em The conferences given at Brown 
University by Solomon Marcus in the years 2008 and 2011}
(in Romanian), Sept. 2025.

 \bibitem{s2} Golopentia, S. {\em When I think at Solomon Marcus}
(in Romanian), Dec.
2025.

\bibitem{s3} Golopentia, S. {\em Reactions to some texts about the
Romanian language of the Acad. Solomon Marcus} (in Romanian),
 August 2013.

\bibitem{s4} Golopentia, S. {Solomon Marcus or about 
the serious joys} (in Romanian), Jan. 2017.

\bibitem{arxive} B.K. Winter, A.T. Lipnicki, 
{\em A braid box}, arXiv: 2604.20884.

\bibitem{bv}
Z.P. Bradshaw and C. Vignat, { \em Learning from Ramanujan:
Elementary Approach to Profound Ideas}. arXiv: 2605.08484.


\bibitem{sm2} Marcus, S. {\em Transcendence, as a universal paradigm} 
(in Romanian), Convorbiri Literare {\bf 2014}, February, 17-28.

\bibitem{sm3} Marcus, S. {\em Transcendence, as a universal paradigm},
BALANCE, A Club of Rome Magazine { 2015}, no.1, March, 50--70.

\bibitem{e} Petrie, B.J. {\em Leonhard Euler's 
use and understanding of mathematical transcendence}, 
Historia Mathematica { 2012}, 39, 280-291.

\bibitem{b2} Nicolescu, B. { \em Transdisciplinarity - 
 past, present and future}, in Moving Worldviews - 
Reshaping sciences, policies and practices for 
endogenous sustainable development { 2006}, COMPAS Editions, Holland, edited by
Bertus Haverkort and 
Coen Reijntjes, 142-166.

\bibitem{bw} Nicolescu, B. { \em The unexplected way to holiness: Simone Weil
(II)} 
 (in Romanian), Convorbiri Literare { 2017}, October, No. 10(262),26--28.

\bibitem{fn} Nichita, F.F. {\em On Models for Transdisciplinarity},
Transdisciplinary Journal of 
Engineering and  Science { 2011}, Vol. 2011, 42-46.

\bibitem{tc} B.E. Tonn and C.J. Hill {\em Capturing the
expanding research areas of the future of humanity
within the field of future studies: The case for a
transcendental future subdiscipline}, Futures 174 (2025) 103692.

\bibitem{racks} Racks and Quandles, Wikipedia, 2026.

\bibitem{ybs} T. Brzezi\'nski and F. F. Nichita,
{\em Yang--Baxter Systems and Entwining Structures},
Comm. Algebra { 33} (2005) 1083-1093.

\bibitem{tt} 
Terence Tao, {\em Mathematics in the age of AI},
2608.16753.

\bibitem{mp} Marta Petreu, { \em Filozofia lui Blaga}, 
Polirom, 2024.

\bibitem{cv} Calin Vlasie, { \em Viitorul literaturii.
Spre o noua poetica a creatiei in era post-algoritmica (I)}, 
Ramuri, 8 / 2026.

\end{thebibliography}
\end{document}